\documentclass[onecolumn]{elsarticle}
\input{amssym}

\begin{document}
\title{Group Analysis via Weak Symmetries For \\ Benjamin-Bona-Mahony Equation }
\author[MN]{M. Nadjafikhah}\ead{m\_nadjafikhah@iust.ac.ir}
\author[]{F. Ahangari}\ead{fa\_ahangari@iust.ac.ir}
\author[]{S. Dodangeh}\ead{s\_dodangeh@mathdep.iust.ac.ir}
\address[MN]{Corresponding author:
School of Mathematics, Iran University of Science
and Technology, Narmak, Tehran 1684613114, Iran.}
\renewcommand{\sectionmark}[1]{}
\begin{abstract}
In this paper, weak symmetries of  the Benjamin-Bona-Mahony (BBM)
equation  have been investigated. Indeed, this method has been
performed by applying the non-classical symmetries of the BBM
equation and the infinitesimal generators of the classical
symmetry algebra of the KdV equation as the starting constraints.
Similarity reduced equations as well as some exact solutions of
the BBM equation are obtained via this method.
\end{abstract}
\maketitle
\section{Introduction}
The Benjamin-Bona-Mahony equation
\begin{eqnarray}
{\rm BBM}\;:\;\;u_t+u_x+uu_x-u_{xxt}=0,\label{eq:1}
\end{eqnarray}
used to model an approximation for surface water waves in a
uniform channel \cite{[9]}. If we note the KdV type equation
\begin{eqnarray}
u_t+u_x+uu_x+u_{xxx}=0,\label{eq:2}
\end{eqnarray}
then we find out the likeness between these equations. Indeed,
this similarity is not stochastic. Both of them used to model the
waves appear in liquids, compressible fluids, cold plasma and
enharmonic crystals which are of surface, hydro-magnetics,
acoustic-gravity and acoustic types, respectively. The interesting
point is that the main difference between equations (\ref{eq:1})
and (\ref{eq:2}) occurs in the case of short waves (Find more
information in \cite{[9],[10]}).

$~~~$ The physical applications and mathematical properties of the
BBM equation (\ref{eq:1}) have been motivated many investigations
such as obtaining the exact solutions via finite difference
discrete process, global attractor and etc.

$~~~$ In this paper, we find the similarity reduced ODEs as well
as resulted similarity solutions of this equation via weak
symmetry implementation. Indeed,  the organization of the present
paper is as follows: Some historical information on the  weak
symmetry method are given in section 2. In section 3, we follow
\cite{[7]} in order to describe the theory of weak symmetries.
Section 4 is devoted to performing this new class of symmetry
methods using the invariant surface condition of the BBM equation
(which is indeed the non-classical symmetry method) and
infinitesimal generators of the classical symmetry algebra of the
KdV  type equation as the starting points in the weak symmetry
method implementation. Finally, we have compared our results with
those related papers using the classical symmetry method in order
to clarify the advantages and disadvantages of the both strong and
weak symmetry methods.
\section{Background}
$~~~$Symmetry methods for differential equations, was originally
developed by {\sc S.~Lie} \cite{[8]}. These methods without any
doubt are very useful and algorithmic  for analyzing and solving
linear and non-linear differential equations. Classification of
differential equations as well as linearization of them are some
other important applications of the symmetry transformation
approach. First {\sc G.W.~Bluman} and {\sc J.D.~Cole} introduced
the notion of the non-classical symmetry group of differential
equations specially for the heat equation in 1969 (Find more
information in \cite{[6]}). For the non-classical method, we seek
the invariance of both the original equation and its invariant
surface condition, exactly this constraint (i.e invariance surface
condition) causes the non-classical solutions which are more
general than the classical ones. There are various implementations
for performing the non-classical symmetry method, for example,
using the compatibility condition has been suggested by {\sc
G.~Cai} and {\sc X.~Ling} \cite{[12]}.

$~~~$ First the weak symmetries have been introduced by {\sc
P.J.~Olver,} and {\sc P.~Rosenau} in 1986 as a generalization of
the non-classical symmetries with motivation of finding every
solutions of the given system. In principle, not only the
invariant solutions corresponding to arbitrary transformation
groups can be found by the reduction method, but also every
possible solution of the system can be found by using some
transformation groups. In other words, there are no conditions
that need to be placed on the transformation group in order to
apply the basic reduction procedure (Find more information in
\cite{[7]}). In the next section, we have an attempt to explain
the notation and implementation of the weak symmetry method by
considering the BBM equation as an example in order to prepare an
appropriate setting.

\section{On the weak symmetry method}
$~~~$Symmetry groups of a system of partial differential equations
can be defined in two types (see \cite{[7]}). %
\paragraph{\bf Definition} {\em
Let $\Delta$ be a system of partial differential equations. A
strong symmetry group of $\Delta$ is a group of transformations
$G$ on the space of independent and dependent variables which has
the following two properties:
\begin{itemize}
\item[a)] The elements
of $G$ transform solutions of the system to other solutions of the
system.
\item[b)] The $G-$invariant solutions of the system are found from a
reduced system of differential equations involving a fewer number
of independent variables than the original system $\Delta$.
\end{itemize}}
\paragraph{\bf Definition} {\em
A weak symmetry group of the system $\Delta$ is a group of
transformations which satisfies the reduction property (b), but no
longer transforms solutions to solutions.}

\medskip $~~~$Indeed, there are several transformation groups which
don't transform solutions of given equations again to solutions,
but their differential invariants enable us to reduce them. In
continuation we would illustrate the procedure of performing this
method. For this purpose, first consider an arbitrary
one-parameter transformation group, then substitute its related
differential invariants and their derivatives into the  original
equation, finally, you will encounter with three different
possible cases which in continuation have been illustrated for the
(BBM) equation using an appropriate one-parameter transformation
group.
\subsection{Reduced equation has  no  parametric variables} Consider the
one-parameter group $$(x,t,u)\mapsto(x+\lambda,t+\lambda,u),$$ So,
we have the characteristic equation $dx=dt=du/0$. By substituting
the resulted differential invariants i.e. $r=x-t$ and $w=u$, into
equation (\ref{eq:1}), we have $w_{rrr}+ww_r=0$. As we see,  this
equation has no parameter variable.
\subsection{Reduced equation isn't incompatible and has
  parametric variables} Consider the one-parameter group
 $$(x,t,u)\mapsto (\lambda x,t,\lambda u).$$ So, the characteristic equation is
$dx/x=dt/0=du/u$. By substituting  the resulted differential
invariants  $r=t$ and $w=u/x$, in equation (\ref{eq:1}), we have
$x(w^2+w_r)+w=0$, where $w=0$ is it's solution and this equation
has $x$ as the parametric variable.
\subsection{Reduced equation is incompatible and has parametric variables} Consider the
one-parameter group $$(x,t,u)\mapsto (x+2\lambda
t+\lambda^2,t+\lambda,u+8\lambda t+4\lambda^2).$$ By substituting
the resulted differential invariants i.e. $r=x-t^2$ and
$w=u-4t^2$, in equation (\ref{eq:1}), we have $
ww_r+w_r+(8-2w_{rrr}-2w_r)t+4w_rt^2$, where this equation has $t$
as the parametric variable and it is incompatible. Indeed, from
the coefficient of $t^2$ we have $w_r=0$ and from the coefficient
of $t$ we have $w_{rrr}+w_r=4$, this means that these equations
are incompatible.
%
\section{Implementation of the weak symmetry method for the BBM equation}
Since, the weak symmetry method is based on conjecture, so here,
the several ideas of performing  this method as well as some of
its aspects are presented.
\subsection{Non-classical symmetries of the BBM equation}
There are several implementations to find the non-classical
symmetries. Here, we follow the procedure presented by {\sc G.
Cai} et al. which they obtained the non-classical symmetries of
the Burgers-Fisher equation based on the compatibility conditions
\cite{[11]}.

Consider the following one-parameter group:
\begin{eqnarray}{\label{eq:7}}
\tilde{x}&=&x+\varepsilon\xi(x,t,u)+O(\varepsilon^2),\nonumber\\
\tilde{t}&=&t+\varepsilon\eta(x,t,u)+O(\varepsilon^2),\\
\tilde{u}&=&u+\varepsilon\varphi(x,t,u)+O(\varepsilon^2),\nonumber\label{eq:3}
\end{eqnarray}
Assume that the equation
$\Delta_1(x,u^{(n)}):=\mbox{eq}(\ref{eq:1})$ is invariant under
the transformation group (\ref{eq:7}) with the following invariant
surface condition:
\begin{eqnarray}{\label{eq:8}}
\Delta_2(x,u^{(n)}):=\eta u_t+\xi u_x-\varphi=0\label{eq:4}
\end{eqnarray}
This means that ${X^{(4)}\Delta_1}|_{\Delta_1=0,\Delta_2=0}=0$,
where
$$X=\xi(x,t,u)\partial_x+\eta(x,t,u)\partial_t+\varphi(x,t,u)\partial_u,$$
is the infinitesimal generator of (\ref{eq:3}), and
$$X^{(4)}=X+\varphi^x\partial_{u_x}+...+\varphi^{tttt}\partial_{u_{tttt}},$$
is the fourth prolongation of $X$, with the coefficients defined
as $\varphi^{J}=D_JQ+\xi u_{Jx}+\eta u_{Jt}$, where $Q=\varphi-\xi
u_x-\eta u_t$ is the Lie characteristic and
$D_J=\sum_{i=0}u_{Ji}\,\partial_{ u_J}$ is the total derivative
w.r.t. $J$ (Find more information in \cite{[1],[2]})

Without loss of generality in condition (\ref{eq:8}), two cases
$\eta=0$ and $\eta=1$ must be considered.

\textbf{Case I $\eta=1$:} In this case we have $u_t=\varphi-\xi
u_x$. Substituting this expression in (\ref{eq:1}) we have
$D_t(\varphi-\xi u_x)=D_t(u_{xxt}-u_x-uu_x)$, where $D_t$ is total
derivative w.r.t. $t$. By substituting $\xi u_{xx}$ in both sides
of above, we find
\begin{eqnarray}
\varphi^t&=&u_{xxtt}-u_{xt}-u_tu_x-uu_{xt}+\xi u_{xx}\nonumber\\
&=&D_{xxt}(u_{t})-(u+1)D_x(u_t)+(\xi u_x-\varphi)u_x+\xi u_{xx},\nonumber\\
&=&D_{xxt}(\varphi-\xi u_x)-(u+1)D_x(\varphi-\xi u_x)+(\xi u_x-\varphi)u_x+\xi u_{xx},\\
&=&\varphi^{xxt}-\xi u_{xxxt}-(u+1)\varphi^x+(u+1)\xi u_{xx}+(\xi
u_x-\varphi)u_x+\xi u_{xx},\nonumber
\end{eqnarray}
By virtue of $D_x(u_t)=D_x(u_{xxt}-u_x-uu_x)$, we have
$u_{xt}=u_{xxxt}-u_{xx}-uu_{xx}-u_x^2$. Finally, we find the
following governing equation:
\begin{eqnarray}
\varphi^t=\varphi^{xxt}-(u+1)\varphi^x-\varphi u_x,\label{eq:5}
\end{eqnarray}
where $\varphi^t=D_t(\varphi-\xi u_x)+\xi u_{xt}$,
$\varphi^x=D_x(\varphi-\xi u_x)+\xi u_{xx}$, and
$$\varphi^{xxt}=D_{xxt}(\varphi-\xi u_x)+\xi u_{xxxt}.$$ By
substituting the coefficient functions $\varphi^t, \varphi^x,
\varphi^{xxt}$ into invariance condition (\ref{eq:5}), we are left
with a polynomial equation involving the various derivatives of
$u(x,t)$ whose coefficients are certain derivatives of $\xi$ and
$\varphi$. Since, $\xi$ and $\varphi$ depend only on $x$, $t$, $u$
we can equate the individual coefficients to zero, leading to
these complete set of determining equations:
$\xi_x=\xi_t=\xi_u=0$, $\varphi=0$. So, we have $\xi=c_1$,
$\varphi=0$. So,  we find the infinitesimal generators of the
non-classical symmetries using the above results as follows, when
$c_1=1$, we have $\sigma_1=u_x+u_t$, and for $c_1\neq0$ the
symmetries are $\sigma_2=u_x$, $\sigma_3=u_t$. As a result we can
state the following proposition:
\paragraph{\bf Proposition} {\em The
non-classical symmetries of the BBM equation in the case of
$\eta=1$, spanned by
\begin{eqnarray}
\sigma_1=u_x+u_t,\qquad\sigma_2=u_x,\qquad \sigma_3=u_t.
\end{eqnarray}}

\medskip As a result of above proposition  we have the following
group-invariant solutions:
\begin{itemize}
\item[1)] For $\mathbf{\sigma_1}=u_x+u_t$, substituting it into
$\sigma_1(u)$ we find $u=F(x-t)$,  where $F$ must  satisfy in:
$FF'-F'''=0$
\item[2)] For $\mathbf{\sigma_2}=u_x$, substituting it into
$\sigma_2(u)=0$ we find $u=F(t)$  for an arbitrary $F$, so from
equation (\ref{eq:1}) we obtain: $u=0$.
\item[3)] For $\mathbf{\sigma_3}=u_t$, substituting it into
$\sigma_3(u)=0$ we find $u=F(x)$, where from equation (\ref{eq:1})
$F$ satisfies this equation: $F'+FF'+F'''=0$.
\end{itemize}
\textbf{Case II $\eta=0$:} In this case,  without lose of
generality we can let $\xi=1$, so we have: $u_x=\varphi$. Using
this we can deduce $A(x,t,u)=\varphi_{xt}-\varphi-u\varphi$.
Subsisting this in the determining equation
$A\varphi_u+\varphi_t-A_u\varphi-A_x=0$, we obtain:
\begin{eqnarray}
\varphi_{xt}\varphi_u-2\varphi\varphi_u-u\varphi\varphi_u+\varphi_t=
\varphi_{xtu}\varphi+u\varphi_u\varphi+2\varphi^2+\varphi_{xxt}+\varphi_x+u\varphi_x.
\end{eqnarray}
By assuming $\varphi=\varphi(x,t)$ above equation changes into
$$\varphi_t-2\varphi^2-\varphi_{xxt}-\varphi_x-u\varphi_x=0.$$ So we
have: $\varphi=1/(c-2x)$. As a result, we deduce that
$u(x,t)=x/(c-2t)+F(t)$ (where $F$ is an arbitrary function) is a
solution of (\ref{eq:1}).
\subsection{Using the classical symmetries of KdV  type equation (\ref{eq:2})}
Since the appearance forms of equation (\ref{eq:1}) and
(\ref{eq:2}) are similar, we want to try our chance in order to
obtain new similarity reduced ODEs for BBM  equation through
infinitesimal generators of the classical symmetries (CS) of KdV
 type equation as the starting constraint. For the classical symmetries
of the KdV  type equation using Lie classical symmetry we have the
next theorem (Since, the proof is computational, to keep scope we
don't present it here. Find more information in \cite{[1],[2]}).
\paragraph{\bf Theorem} {\em If we consider
$X=\xi(x,t,u)\partial_x+\eta(x,t,u)\partial_t+\eta(x,t,u)\partial_u$
as the infinitesimal generator of the classical symmetry group of
the KdV  type equation (\ref{eq:2}),  then we have
\begin{eqnarray}
\eta=c_1t+c_2,\qquad\xi=\frac{1}{3}c_1(x+2t)+c_3t+c_4,\qquad
\varphi=-\frac{2}{3} c_1u+c_3,\label{eq:6}
\end{eqnarray}
where $c_1$, $c_2$, $c_3$ and $c_4$ are arbitrary constants.}

\medskip Hence the next corollary could be stated:
\paragraph{\bf Corollary} {\em
 The classical symmetries of equation (\ref{eq:2}) i.e. KdV  type equation,
spanned by:
\begin{eqnarray}
X_1=(x+2t)\partial_x+3t\partial_t-2u\partial_u,\;\;
X_2=\partial_t,\;\; X_3=t\partial_x+\partial_u,\;\;
X_4=\partial_x.
\end{eqnarray}}

$~~~$So, we can consider any linear combinations of given vector
fields in the above corollary as the starting constraint of the
weak symmetry method. In continuation, we will illustrate the weak
symmetry method using some linear combination of $X_1$, $X_2$,
$X_3$ and $X_4$ as the starting point.
\paragraph{\bf Example} Consider the one-parameter
transformation group with the infinitesimal generator
$X_2+X_3=t\partial_x+\partial_t+\partial_u$. The characteristic
equation is $dx/t=dt=du$. So, we find the differential invariants
as $r=t^2-2x$, $w=u/t$. By substituting these new variables in the
original equation (\ref{eq:1}) we deduce
$(2w_r-2ww_r-w_{rrr})t^2=2w_rt+4w_{rr}+w$, where $t$ can be
considered as the differential parameter. Note that solving the
above ODE doesn't give new solution.
\paragraph{\bf Example}
Consider the one-parameter transformation group with the
infinitesimal generator $X_3=t\partial_x+\partial_u$. The
characteristic equation is $dx/t=dt/0=du$. So, we can obtain the
differential invariants as $r=t$, $w=u-x/t$. By substituting these
new variables in the original equation (\ref{eq:1}) we find:
$rw_r+w-1=0,$, solving this reduced equation we obtain
$w=r/(r+c)$. So we can find $u=(tx+x^2+cx)/(t(x+c))$ as the
solution of equation (\ref{eq:1}).
\subsection{Some other suggestions}
Some other ideas may be useful to reach other solutions of the BBM
equation. For example, non-classical potential symmetry method or
using classical and non-classical symmetries of other equations
which have the similar forms as the BBM equation. Meanwhile,
Physical knowledge of the model framework can be so effective in
order to reach favorite solutions via weak symmetries. For example
if you know your desired solution may be invariant under some
scale of specific variables then the weak symmetry method can be
started with an appropriate scaling transformation. Since the main
goal of this paper was introducing weak symmetry method for BBM
equation, we lay away  performing of above approaches.
\section{More discussions}
Now, we want to compare our results with other related papers.
Paper \cite{[5]} is concentrated on the classical symmetries and
optimal Lie system of the BBM equation. Comparing with \cite{[5]},
we deduce that in this paper by applying the weak symmetry method
we have obtained more similarity solutions  and other useful
suggestions are presented in order to reach more other
solutions.\newline
$~~~$ Taking into account the sections 2 and 3 of \cite{[5]}, the
next theorem can be resulted (Find more information in
(\cite{[1]}, Chapter 3).
\paragraph{\bf Theorem} {\em If $u=f(x,t)$ is solution of the BBM equation (\ref{eq:1}), so are
the functions
\begin{eqnarray}\nonumber
u=f(x-\varepsilon,t),\quad
u=f(x-\alpha\varepsilon,t-\varepsilon),\quad
u=e^{(u+1)\varepsilon}f(x-\alpha\varepsilon,e^{-\varepsilon t}t),
\end{eqnarray}
where $\varepsilon\ll1$ and $\alpha$ are arbitrary constants.}

\medskip Indeed, above theorem characterizes the invariant solutions
of the BBM equation, for instance if $u=c$ is a solution of
equation (\ref{eq:1}), then  from this theorem we obtain
$u=ce^{\varepsilon(u+1)}$ as a solution of the BBM equation. For
another example, if we consider the solution
$u=(tx+x^2+cx)/(t(x+c))$ of equation (\ref{eq:1}), from this
theorem we deduce  that
$$u=\displaystyle{\frac{(t-\varepsilon)(x-\alpha
\varepsilon)+(x-\alpha\varepsilon)^2+c(x-\alpha\varepsilon)}
{(t-\varepsilon)(x-\alpha\varepsilon+c)}},$$ (where
$\varepsilon\ll1$ and $\alpha$, $c$ are arbitrary constants), is
again a solution of BBM equation. By using such approach, we are
enable to obtain  more new solutions for the BBM equation.
\section*{Conclusions}
$~~~$In this paper,  we have presented a comprehensive explanation
of the weak symmetry method as the generalization of the classical
Lie symmetry method. Indeed, we have performed the weak symmetry
method for the BBM equation which has been fulfilled by applying
the non-classical symmetries of the BBM equation and using the
classical symmetries of the KdV  type equation as the starting
constraints. Also, the similarity reduced equations as well as
some exact solutions of the BBM equation are obtained via this
method. Finally, we have compared our results with papers using
the classical symmetry method. Other suggestions for finding new
exact solutions are also presented.

%

\end{document}